\documentclass[12pt]{amsart}
\usepackage{amsthm,amsmath,amssymb,amscd,graphics,enumerate}

{
   \newtheorem{theorem}[subsubsection]{Theorem}
   \newtheorem{proposition}[subsubsection]{Proposition}     
   \newtheorem{lemma}[subsubsection]{Lemma}

   \newtheorem{corollary}[subsubsection]{Corollary}

}
{\theoremstyle{definition}

}
\newcommand{\bbS}{{\mathbb{S}}}

\newcommand{\PP}{{\mathbb{P}}}

\newcommand{\bK}{{\mathbf{K}}}
\newcommand{\bM}{{\mathbf{M}}}

\newcommand{\cC}{{\mathcal C}}
\newcommand{\cD}{{\mathcal D}}
\newcommand{\cE}{{\mathcal E}}

\newcommand{\cK}{{\mathcal K}}

\newcommand{\cM}{{\mathcal M}}

\newcommand{\cO}{{\mathcal O}}

\newcommand{\cS}{{\mathcal S}}

\newcommand{\cX}{{\mathcal X}}

\newcommand{\Spec}{\operatorname{Spec}}

\newcommand{\st}{{\operatorname{st}}}

\newcommand{\das}{\dashrightarrow}

\newcommand{\can}{{\operatorname{can}}}
\newcommand{\bal}{{\operatorname{bal}}}

\newcommand{\red}{{\operatorname{red}}}
\newcommand{\norm}{^{\operatorname{norm}}}

\newcommand{\double}{\genfrac..{0pt}1
{\raise -1pt\hbox{$\scriptstyle\longrightarrow$}}{\raise 3pt\hbox
{$\scriptstyle\longrightarrow$}}} 
\newcommand{\setmin}{\,\protect%
\begin{picture}(8,3.5)\qbezier(1,3.5)(4,2.)(7,.5)\end{picture}\,}

\renewcommand{\setminus}{\setmin}

\newcommand{\KO}[4]{{\cK_{#1,#2}(#3,#4)}}
\newcommand{\KOB}[4]{{\cK^{\bal}_{#1,#2}(#3,#4)}}
\newcommand{\ko}[4]{{\bK_{#1,#2}(#3,#4)}}

\newcommand{\TSM}{\KO{g}{n}{\cM}{d}}
\newcommand{\TSMM}{\KO{g}{\{n\}}{\cM}{d}}
\newcommand{\tsm}{\ko{g}{n}{\cM}{d}}

\newcommand{\TSMB}{\KOB{g}{n}{\cM}{d}}

\def\tototi{\mathbin{\mathop{\otimes}\limits^{\raise-1pt\hbox
{$\scriptscriptstyle {\rm L}$}}}}

\def\indlim{\mathop{\vrule width0pt height7pt depth
4pt\smash{\lim\limits_{\raise 1pt\hbox to 14.5pt
{\rightarrowfill}}}}}
\def\projlim{\mathop{\vrule width0pt height7pt depth
4pt\smash{\lim\limits_{\raise 1pt\hbox to 14.5pt
{\leftarrowfill}}}}}

\begin{document}

\title[Plurifibered varieties]{Canonical models and stable reduction for
plurifibered varieties} 
\author[D. Abramovich]{Dan Abramovich}
\thanks{Research partially supported by NSF grant DMS-0070970, a
Forscheimer Fellowship and Landau Center Fellowship}  
\address{Department of Mathematics\\ Boston University\\ 111 Cummington
         Street\\ Boston, MA 02215\\ U.S.A.} 
\email{abrmovic@math.bu.edu}\maketitle

\section{Introduction}
All schemes considered here are over some noetherian base scheme $\bbS$ of pure   
characteristic 0.
\subsection{Prologue: semistable reduction over a base of higher dimension}

One result we prove here is 
the
following weak variant of the main theorem of \cite{A-K}:

\begin{theorem}\label{Th:semistable} Let $X \to B$ be a morphism of
projective varieties, with 
geometrically integral generic fiber. Then there is an alteration $B_1
\to B$, a modification of the main component 
$Y\to X\times_B B_1$, and a divisor $D\subset Y$  such that 
\begin{enumerate}
\item both $Y
\to B_1$ and $D\to B_1$ are flat, 
\item the singularities of $(Y, D) $ are log-canonical, and
\item  for every
$b\in B_1$, the fiber $(Y_b, D_b)$ has at most semi-log-canonical
singularities. 
\end{enumerate}
\end{theorem}

This result is weaker that \cite{A-K}, Theorem 1, since the
singularities of $(Y, D) $ are only log-canonical, and we do not claim
that the singularities of $Y$ itself are Gorenstein, or toroidal, or
canonical. However, as we will see, the procedure that leads to $(Y,D)$ is more
canonical than that of \cite{A-K}, since the only data required is a
rational plurifibration, and  once this is chosen everything is
canonically defined. The main idea is here inspired by
A.J.\,de Jong's work on alterations \cite{DeJong}.

The contents of this note are centered around stable models and stable
reduction for plurifibrations. This  is an expansion of remarks which apeared
in 
\cite{stablemaps} and \cite{modfam}, on natural extensions of the
results of \cite{fibsurf}. Related ideas, but with quite a
different bent, were independently discovered by Mochizuki
\cite{Mochizuki}.

\subsection{The minimal model program}
One could say that the basic goal of the minimal model program (MMP)
is to find, in the birational equivalence class of each variety $X$ of
general type, a canonical model 
$$X \das X^\can,$$ where $X^\can$ has canonical singularities and ample
canonical divisor.  A more ambitious goal proposed by V. Alexeev is that
of stable maps: given a triple
$$( X, D, f: X \to \PP),$$
where $X$ is, say, a smooth projective variety, $D$ a normal-crossings
divisors, and $f:  X \to \PP$ is a morphism to a fixed projective
space, with the ``pre-stability assumption'' $$L_{( X, D, f: X \to
\PP),n}:=\ K_X(D) \otimes
f^*\cO_\PP(n)\ \ \mbox{ is a big line bundle for }\ n >> 0, $$
one expects to find a suitable birational representative
$$( X^\st, D^\st, f^{\st}: X \to \PP)$$ which is a {\em stable map},
that is, $( X^\st, D^\st)$ log-canonical and $$K_{X^\st}(D^\st) \otimes
(f^\st)^*\cO_\PP(n)\ \ \mbox{ is an {\bf ample} line bundle for } \ n>>0. $$ 

Following Deligne-Mumford \cite{Deligne-Mumford},
Alexeev\cite{Alexeev:stable-maps}, 
Koll\'ar, Shepherd-Barron \cite{Kollar,Kollar-Shepherd-Barron} and
Viehweg\cite{Viehweg}, one also considers moduli. One defines certain
degenerate 
versions of stable $( X, D, f: X \to \PP)$  above, where $(X, D)$
becomes {\em semi-log-canonical}. Then one conjectures that, fixing
suitable numerical invariants, stable
maps admit a proper moduli stack having a projective coarse moduli
space. It should be pointed out that, following Koll\'ar \cite{Kollar},
Alexeev\cite{Alexeev:stable-maps} and 
Karu \cite{Karu}, the projectivity of moduli in a sense follows from a strong 
enough version of the MMP in one dimension higher.

\subsection{Tweaking with the data}
The stable models above are supposed to arise by a prescribed
procedure: $X^\st$ is the projective spectrum of the algebra of
sections  $$\otimes_{i\geq 0}L_{( X, D, f: X \to
\PP),n}^{\otimes i},$$ and the divisor $D^\st$ and morphism $f^\st$
are induced. Without this requirement, one can easily cheat, i.e by
throwing in a sufficiently  ample divisor $D^\st$, or by replacing
$f$ by an embedding. Moreover, such ``cheat solutions'' are far from
canonical, are not birational in nature, and certainly are not
sufficient in 
finding a good compactification of moduli. Nevertheless, the idea of
temporarily ``tweaking'' with the given data has been quite fruitful:
methods of changing the divisor $D$ a little bit are fundamental in the
work of an illustrious list of authors, specifically Shokurov; and the
addition of a map $f:X\to Z$ 
(specifically a small contraction of an extremal ray)
puts the study of flips right into this framework.

In this note, we propose a situation where the MMP has a complete
solution - stable models and complete moduli included - which is cheap
given existing theory (though the existing theory is far from
cheap). It is still a ``cheat'' solution, but in a milder sense: the
additional data we add to the given variety is {\em birational} in nature,
and the final maps and divisors arise {\em canonically} from the
data, no further choices needed. As with other ``cheat'' solutions,
one hopes that this can serve as a stepping stone for further
undferstanding of the MMP. A situation where this can be
done is studied in the paper \cite{Lanave}.

\subsection{The setup}
We are given a function field $K$ of transcendence
degree $d$ over a field $k$ of characteristic 0. A {\em rational
plurifibration} $(K_i)$ is a sequence 
$$ K = K_0 \supset K_1\supset \cdots \supset K_d = k$$
where each $K_i \supset K_{i+1}$ is an extension
of transcendence degree 1 such that $K_{i+1}$ is algebraically closed
in $K_i$. We denote by $C_i \to \Spec K_{i+1}$ the smooth projective
curve associated to the function field extension $K_i \supset K_{i+1}$. 

Suppose in addition we are given a 
 projective variety $\bM_0$. A {\em rational
plurifibered map} with target $\bM_0$ is a  rational
 plurifibration $(K_i)$ and a morphism 
$\Spec K_0\to \bM_0$. 

A slightly more general notion is the following: suppose $\cM_0$ is a
proper Deligne--Mumford stack having $\bM_0$ as its coarse moduli
space. Then   a {\em rational
plurifibered map} with target $\cM_0$  is a marked rational
plurifibration $(K_i)$ 
and an object
$\cM_0(\Spec K_0)$ (equivalently, a morphism $\Spec K_0\to \cM_0$).

We need to impose a pre-stability condition. The most general
situation allowed is a bit subtle to describe - essentially it is
tautologically  the situation where the main result works. In order to
be more specific we will use  a sufficient condition for
pre-stability: we say  that $(K_i,\Spec K_0\to \cM_0 )$ is {\em
sufficiently good}
if one of the following conditions holds:
\begin{enumerate}
\item[\bf G1.] the morphism $\Spec K_0\to \cM_0$ is finite to its image, or
\item[\bf G2.] the genus $g(C_i)>1$, for $i=0,\ldots,d-2$, and
either $g(C_{d-1})>1$ or $C_{d-2} \to \Spec K_{d-1}$ is nonisotrivial, or
\item[\bf G3.] the function field $K_0$ is of general type over $K_d$.
\end{enumerate}

There are some natural extensions, where
$C_i$ are allowed to have a-priori markings and where more maps $\Spec
K_i \to \cM_i$ are a-priori given. We delay facing these to a later version.

\subsection{Twisted stable maps} Our main tool is the moduli stack of
{\em twisted stable maps}. Recall that a {\em twisted curve} is a
marked Deligne--Mumford stack of dimension 1 with simple local structure
(\cite{stablemaps}, Definition 4.1.2).  A {\em twisted stable map} $f: \cC
\to \cM$ is a representable morphism from a twisted curve $\cC$ to a
Deligne--Mumford stack $\cM$ such that the resulting map $C \to \bM$
of course moduli spaces is stable (\cite{stablemaps}, Definition
4.3.1). It is proven in \cite{stablemaps}, Theorem 1.4.1 that
families of twisted stable with fixed genus $g$, number of markings $n$, and image
class $\beta$  form a Deligne--Mumford stack $\TSM$ having a
projective coarse moduli space $\tsm$. In our discussions below we
will restrict attension to the stack of {\em balanced} twisted stable
maps $\TSMB$. From now on we will assume this and suppress the
superscript ``$\bal$''. In addition, we will ignore the ordering of the markings,
thus consider the stack $$\TSMM := [\TSM/\cS_n],$$ where the symmetric group
$\cS_n$ acts by permuting the numbers of the markings. Thus the stack $\TSMM$
parametrizes maps where the twisted curves $\cC$ are marked by
$\Sigma^\cC$, which is a reduced Cartier divisor which is $n$-sheeted finite
\'etale over the base. In most of what follows the
marking $\Sigma^\cC$ will all be twisted, i.e. the inertia group at
any point of $\Sigma^\cC$ is nontrivial.

\subsection{The outcome: twisted version}
We now describe the type of objects that arise as the stable models in
a certain sense of a rational plurifibered map.

By a {\em twisted stable plurifibered map} over a scheme $S$ we mean  a gadget
$$ (\cX_i \to \cX_{i+1},\Sigma_i, \cX_i\to \cM_i)$$
where 
\begin{enumerate}
\item $\cX_i$ are Deligne--Mumford stacks, $i=0,\ldots,d-1$, and $\cX_d
= S$
\item $ (\pi_1:\cX_i \to \cX_{i+1},\Sigma_i\subset\cX_i, f_i:\cX_i\to
\cM_i)$ is a family 
of twisted
stable maps of certain genus $g_i$, unordered twisted marking of
degree $n_i$, and 
image class $\beta_i$ for each $i=0,\ldots,d-1$, and
\item $\cM_{i+1} = \cK_{g_i, \{n_i\}, \beta_i}$ for $i=0,\ldots,d-2$.
\end{enumerate}

We remark that, almost by definition of twisted stable maps, each $\cX_i$ has a
representable dense open substack, and the coarse moduli space $X_i$
is projective over $S$. 

We use the notation $q_i: \cX_i \to X_i$
for the canonical map.
Also, for each $0\leq i\leq j\leq d$ we denote by $\pi_{i,j}: \cX_i \to
\cX_j$ the composite $\pi_{j-1}\circ\cdots \circ \pi_i$.

\subsection{Main results: twisted version}

The existence of stable models is given by the following theorem:

\begin{theorem}\label{Th:twisted} Let $(K_i, f: \Spec K_0 \to \cM_0)$
be a sufficiently good rational 
plurifibered map. Then
there exists a canonical  twisted stable 
plurifibered map $ (\pi_1:\cX_i \to \cX_{i+1},\Sigma_i\subset\cX_i,
f_i:\cX_i\to 
\cM_i)$ over $S=\Spec K_d$,  where
$\cX_i$ is irreducible with function field $K_i$, with $\pi_i$
corresponding to $K_i\supset K_{i+1}$.
\end{theorem}

Completeness of moduli is immediately given:
\begin{theorem} The category of families of twisted stable
plurifibered maps with fixed data $\cM_0, g_i, n_i, \beta_i$ is a
proper Deligne--Mumford stack having a projective coarse moduli space.
\end{theorem}
{\em Proof.} This category is equivalent to
$\cK_{g_{d-1},\{n_{d-1}\}}(\cM_{d-1},\beta_{d-1})$, for which the
claim follows from \cite{stablemaps}, Theorem 1.4.1. \qed

\subsection{Main results: coarse version} Let $ (\pi_1:\cX_i \to
\cX_{i+1},\Sigma_i\subset\cX_i, f_i:\cX_i\to 
\cM_i)$ be a twisted stable plurifibered variety.

We define a divisor on $X_0 = X$: for each $i$ we denote by
$\Sigma^{X_i}$ the marking on $X_i$ thought of as the coarse curve of
$\cX_i \to \cX_{i+1}$.  
Let $D\subset X$ be given by
$$D = \sum_{i=0}^{d-1} (\pi_{0,i}^*\Sigma^{X_i})_\red.$$

We define a projective variety $\bM$ and a morphism $f:X \to \bM$: 
let $\bM_i$ be the  coarse moduli space of $\cM_i$ and let 
$$\bM = \prod_{i=0}^{d-1} \bM_i.$$
We define $X \to \bM$ to be the product of the morphisms obtained by
composing $X \to X_i$ with $X_i \to \bM_i$.

\begin{proposition}\label{Prop:alexeev-sm} 
The pair $(X,D)$ is semi-log-canonical, and the triple $(X,D,f:X \to
\bM)$ is a stable map in the sense of Alexeev. 
\end{proposition} 

  We call $(X,D,f:X \to \bM)$ the Alexeev stable map associated to $
  (\pi_1:\cX_i \to  \cX_{i+1},\Sigma_i\subset\cX_i, f_i:\cX_i\to  
\cM_i)$.

If $(K_i, D_i, f: \Spec K_0 \to \cM_0)$ is a sufficiently good
rational plurifibered 
map, and $(\pi_1:\cX_i \to
\cX_{i+1},\Sigma_i\subset\cX_i, f_i:\cX_i\to 
\cM_i)$ the associated twisted stable plurifibered map, then we obtain:

\begin{corollary} Given a sufficiently good  rational plurifibered
map, with $K_0$ of 
general type over $K_d$,  there is a
canonical Alexeev stable map $(X,D,f:X \to
\bM)$ associated to the corresponding twisted stable plurifibered map.
\end{corollary}

\section{Proofs}
\subsection{Proof of Theorem \ref{Th:twisted}} The proof proceeds by
induction on $d$. For $d=0$ there is nothing to prove. Assume the
claim is proven for $d-1$, and denote the stable twisted plurifibered
map over $K_{d-1}$ associated to $$K_0\supset\cdots \supset K_{d-1}$$
by $$(\pi_i^{d-1}:\cX_i^{d-1} \to\cX_{i+1}^{d-1}, \Sigma_i^{d-1},
f_i^{d-1}: \cX_i^{d-1} \to \cM_i).$$

{\sc Construction of target.} Denoting $\cM_{d-1} =
\cK_{g_{d-2},\{n_{d-2}\}}(\cM_{d-2}, \beta_{d-2})$, 
we have by definition a canonical morphism $\Spec K_{d-1} \to
\cM_{d-1}$.

{\sc Construction of twisted curve and map.} Essentially by
\cite{stablemaps}, Lemma 
7.2.6 there exists a 
unique smooth twisted curve $\cX_{d-1}\to \Spec K_d$ having function
field $K_{d-1}$, with a representable morphism $f_{d-1}:\cX_{d-1} \to
\cM_{d-1}$ and all markings twisted. Since that lemma had extra
assumptions, we recall the 
construction: the morphism $\Spec K_{d-1} \to \bM_{d-1} $ extends
uniquely to $C_{d-1} = X_{d-1} \to \bM_{d-1}$, since the latter  is
proper. We have a morphism $$\psi:\Spec K_{d-1} \to
C_{d-1}\times_{\bM_{d-1}}\cM_{d-1}.$$ 
Let $$\cX_{d-1} = (\overline{ \psi(\Spec K_{d-1} )} )\norm,$$
i.e. the normalization of the closure of the image of $\psi$. This is
the required twisted curve. 

{\sc Numerical invariants.} We set $g_{d-1} = g(C_{d-1}); n_{d-1} = $
the degree of the marking on $\cX_{d-1}$, and $\beta_{d-1} = $ the
class of $f_{d-1\,*}[C_{d-1}]$.

{\sc Construction of $\cX_i\to \cM_i$.} We construct $\cX_i$ by descending
recursion, $\cX_{d-1}$ being already given. Assume $\cX_{i+1}\to
\cM_{i+1}$ given, we take $(\cX_i\to \cX_{i+1}, \Sigma_i,\cX_{i}\to
\cM_{i})$ to be the universal object corresponding to $\cM_i =
\cK_{g_i, \{n_i\}}(\cM_i, \beta_i)$. We are also given canonical isomorphisms 
$$(\cX_i\to \cX_{i+1}, \Sigma_i,\cX_{i}\to
\cM_{i})_{\Spec K_{d-1}} \ \ \cong \ \  (\cX_i^{d-1}\to
\cX_{i+1}^{d-1}, \Sigma_i^{d-1},\cX_{i}^{d-1}\to 
\cM_{i}), $$ as required.  

{\sc Stability of map.} We claim that $\cX_{d-1} \to \cM_{d-1}$ is
stable. We need to show: 
\begin{quote}  if $X_{d-1} \to \bM_{d-1}$ is constant,
then $2g_{d-1} - 2 + n_{d-1} > 0$, namely:
\begin{enumerate}
\item[claim 1:] if $g_{d-1} = 1$ then $ n_{d-1} > 0$, and
\item[claim 2:] if $g_{d-1} = 0$ then $ n_{d-1} > 2$.
\end{enumerate}
\end{quote} The claims are obvious under the ``sufficiently good''
assumptions (G1) and (G2). We need 
to prove claims (1),(2) under assumtion (G3) that $K_0$ is of general type
over $K_d$. The claims, as well as the general type assumption, are
geometric, so we may replace $K_d$ by its algebraic closure. We prove
the claim by contradiction. 
We now construct trivializations of $\cX_{d-1} \to \cM_{d-1}$ after
finite \'etale covers.

 Assuming we have  $g_{d-1} = 1$ and $
n_{d-1} = 0$, then $\cX_{d-1} = X_{d-1}$ is a smooth elliptic
curve. Since $X_{d-1} \to \bM_{d-1}$ is constant, the object
corresponding to $X_{d-1} \to \cM_{d-1}$ is isotrivial, and becomes
trivial after a finite \'etale cover $E \to X_{d-1}$, and therefore
$E$ is an elliptic curve. 

 Assume, on the other hand, we have 
  $g_{d-1} = 0$ and $
n_{d-1} <2$. Then $X_{d-1}$ is rational and has two marked
  points. Since we have made the base field algebraically closed, we
may choose coordinates so that the points are $0$ and $\infty$. The
  map $X_{d-1} \setminus \{0, \infty\} \to \cM_{d-1}$ is isotrivial,
  which means that there is a map $E = \PP^1 \to X_{d-1}$ (branched only
  at $0$ and $\infty$) with a lifting $\PP^1 \to \cM_{d-1}$. It
  actually follows that this map can be chosen \'etale over
  $\cX_{d-1}$, but this will not be essential.

In both cases we got a trivialization over a curve $E$ with genus
$\leq 1$. By decending induction and the construction
of $E$ and $X_{i}$ we obtain that $E \times_{X_{d-1}} X_0 \to E$ is  a trivial 
family $E \times Y \to E$ . We have a finite morphism  $E
\times_{X_{d-1}} X_0 \to X_0$, and since the source is covered by
curves $E$ of genus $\leq 1$ we get a contradiction to the assumption that
$X_0$ is of general type. This completes claims (1),(2) and proves that
$\cX_{d-1} \to \cM_{d-1}$ is stable.

Thus the theorem is proven. \qed

\subsection{Proof of Proposition \ref{Prop:alexeev-sm}}

We need to show that $(X,D)$ has semi-log-canonical singularities and
that a certain line bundle is ample. It is convenient to study first
the situation on $\cX=\cX_0$, which has an analogous divisor 
$$\cD = \sum_{i=0}^{d-1} \pi_{0,i}^*\Sigma_i.$$ This divisor is
automatically reduced. 
\begin{lemma} The pair $(\cX, \cD)$ has semi-log-canonical
singularities. Moreover, when $\cX$ is normal, the pair has
log-canonical singularities.
\end{lemma}
{\em Proof.}  We use induction on $d$. Since $\cX$ is a composite of
a fibration by nodal curves it is Gorenstein, in particular Cohen-Macaulay.
Thus to prove the claim it suffices to look at the normalization, as follows:

Let $\cX_0\norm \to \cX_0$ be the
normalization, and let $\cE\subset \cX_0\norm$ be the conductor
divisor. Consider an irreducible component $\cX'
\subset\cX_0\norm$. Denote $\cD' =\cD|_{\cX'}+\cE$.   It suffices to
show that the pair $(\cX',\cD')$ has 
log-canonical singularities.  

By definition $\cX'$ is a composite of
fibrations by twisted curves $(\cX_i' \to \cX_{i+1}', \Sigma_i')$,
each with smooth generic fiber  (In particular $\cX_{d-1}'$ is a smooth
curve). 
By
\cite{fiberproduct}, Corollary 4.5,  $\cX'$ has rational Gorenstein
singularities, in particular canonical.
 We are going to 
show that $(\cX', \frac{m-1}{m} \cD')$ is log-canonical for every positive
integer $m$. Taking the limit as $m \to \infty$ we get that the pair $(\cX',
\cD')$ is log-canonical. 

   So fix an integer $m$. Let $p_0$ be a 
geometric point on  
$\cX_0'$ and $p_i$ its image in   
$\cX_i'$.  Passing to an \'etale neighborhood of $p_0$  we may replace $\cX'$
by a 
scheme, which we still denote by $\cX'$. For each $i$ such that $p_i$ lies on
$\Sigma_i'$ we have that $\cX_i 
\to \cX_{i+1}$ is smooth and $\Sigma_i'$ is a section. Passing to a
neighborhood again, we may assume there exists a cyclic cover $Y_i \to
\cX_i'$ of degree $m$ totally branched along $\Sigma_i'$ and \'etale
ensewhere. Pulling back and taking the compositum, we get a branched cover
$\psi_i:Z_i 
\to \cX_i'$, totally branched along all the components of $\cD_i$ with index
$m$, and each $Z_i \to Z_{i+1}$ is a family of nodal curves. Again by
\cite{fiberproduct}, Corollary 4.5, $Z_0$ has canonical singularities.

Note that $K_{Z_0} = \psi_0^* (K_{\cX'} + \frac{m-1}{m}\cD')$. Since $Z_0$ has
canonical singularities, a result of
Koll\'ar (see \cite{flab}, 20.3 (2)) says that the pair $(\cX', \cD')$ is
log-canonical.

\qed

We continue with the proof of the proposition. The next step is
\begin{lemma} The pair $(X,D)$ has semi-log-canonical
singularities. Moreover, for a positive integer $m$ such that
$O_X(m(K_X+D))$ is Cartier, we have
$p^*(O_X(m(K_X+D))) = (\omega_{\cX}(\cD))^m$.   
\end{lemma} 
{\em Proof.}
Recall that when $\cX$
is a Deligne--Mumford stack with coarse moduli space  $X$, there is an
\'etale covering $\cup X_j \to X$, schemes $V_j$, and finite groups
$\Gamma_j$ such that $\cX\times_X X_j \cong [V_j / \Gamma_j]$. In our
case, we already have $(\cX,\cD)$ semi-log-canonical, so $(V_j,
\cD_{V_j})$ are semi-log-canonical. We also have that $\cD_{V_j}$
contains the divisorial part of the fixed-point-logus of
$\Gamma_j$. Denote by $D_j$ the image of $\cD_{V_j}$ in $V_j/\Gamma_j$. The
aforementioned result of Koll\'ar (see \cite{flab}, 20.3 (2))  implies that
$(V_j/\Gamma_j,D_j)$ has 
semi-log-canonical singularities, and the lemma follows. \qed 

Let $\cO_\bM(1)$ be an ample line bundle.
The proposition follows once we prove the following lemma:
\begin{lemma} For any large enough integer $n$
the sheaf $\omega_{\cX}(\cD) \otimes f^*\cO_\bM(n)$ is ample.
\end{lemma} 

{\em Proof.} There is an integer $n$ such that
$L_i : = \omega_{\cX_i/\cX_{i+1}}(\Sigma_i) \otimes f_i^*\cO_{\bM_i}(n)$ is
$\pi_i$-ample. It follows from a result of Koll\'ar (\cite{Kollar}, Proposition
4.7)  
that  $L_i$ is nef on $\cX_i/S$: indeed, $L_i^m$ is $\pi_i$-ample and
$\pi_i$-acyclic for large
$m$; Koll\'ar shows that the resulting vector bundle $\pi_{i\,*}L_i^m$ is
 semipositive, therefore the tautological line-bundle on its projectivization
is 
nef, and therefore the restriction of this tautological bundle to 
$\cX_i$, which is  $L_i^m$, is nef.

Let $F \subset \cX$ be closed and irreducible of dimension $l$. We
need to show that $$\deg_F c_1(L)^l >0.$$

We use induction on $d$ (the case $d=0$ being trivial).

If $\dim( \pi_0(F)\subset \cX_1) = l$ then we can use the inductive
assumpltion on $d$ applied to $\cX_1 \to S$ and the fact that $L_0$ is
nef to conclude that $\deg_F c_1(L)^l >0.$ 

Otherwise we have  $\dim( \pi_0(F)) = l-1$. Denote $\pi_0(F) = F_1$.
Write $L = L_0 \otimes \pi_0^* L'$, where $L'$ is the analogous line
bundle for $\cX_1/S$. By assumption $\deg_F \pi_0^*c_1(L')^{l}=0$, but
the inductive assumption implies $\deg_{F_1}c_1(L')^{l-1}>0$.
Now using the projection formula
$$\deg_F c_1(L)^l > \deg_F (c_1(L_0)\cdot \pi_0^*c_1(L')^{l-1}) = \deg_{F_1}(
\pi_{F/F_1\,*}c_1(L_0)\ \cdot\ c_1(L')^{l-1})>0.$$
By the Kleiman-Nakai-Moishezon criterion for ampleness this completes
the Lemma.\qed 
\subsection{Proof of Theorem \ref{Th:semistable}}
Let $X\to B$ be as in the theorem. We set $\bbS = B$ as base scheme,
and $\bM_0 = X$ as target scheme. We set $K_0 = K(X), K_d = K(B)$, and
write $\eta = \Spec K_d$. Now
define $K_i,\ i=1,\ldots,d-1$ by descending recursion: given $K_i$ choose a
general pencil of hypersurfaces $H_t \to\PP^1_{K_i}$ of $X_{K_i}$. We
have $K(H) = K_0$, and we set $K_{i-1} = K( \PP^1_{K_i})$. (The
connectedness theorem implies that $X_{K_i}$ is geometrically integral
for all $i$.)

Since the morphism $$\Spec K_0 \to \bM = X$$ is an embedding, the rational
plurifibered map $(K_i, \Spec K_0 \to \bM_0)$ is sufficiently good. We
therefore have a twisted stable plurifibered map over $\eta$
$$(\cX_i^\eta, \Sigma_i^\eta, f_i^\eta: \cX_i^\eta \to \cM_i).$$

This gives a morphism $$\eta \to \cK_{g_{d-1}, \{n_{d-1}\}}(\cM_{d-1},
\beta_{d-1}).$$ Since the latter is proper, there exists an alteration
$B_1 \to B$ and a lifting 
$$B_1\to \cK_{g_{d-1}, \{n_{d-1}\}}(\cM_{d-1},
\beta_{d-1}),$$ in other words, a family of  stable plurifibered maps
over $B_1$
$$(\cX_i^1, \Sigma_i^1, f_i^1: \cX_i^1 \to \cM_i)$$ 
(extending the pullback of the family over $\eta$ to the generic point
of $B_1$). In particular we have a rational map $\cX_0 \to X$. 

We now have
 an associated Alexeev stable map
$$(Y, D, Y \to \bM),$$ 
where by construction $\bM = \bM_0 \times \bM'$ for some projective scheme
$\bM'$. A-priori we have a rational map $Y \to X$, but the
composition $Y \to \bM \to \bM_0= X$ shows that this is regular.
\qed
\subsection{Acknowledgements} I wish to thank the Hebrew University and its
Lady Davis Foundation and Landau Center for Mathematical Research for support
of a visit during which this note was written. I also thank the organizers of
the HDG program at the Newton Institute for an invitation for a visit during
which this note was completed, and for a number of participants, including
Koll\'ar, Mosta\c{t}a and Reid, who clarified a few murky points about
inversion of adjunction. Finally I thank 
A. Vistoli for many discussions over the years which set the stage for the
main idea here.


\begin{thebibliography}{MMMM}
\bibitem[$\aleph$]{fiberproduct} D. Abramovich, {\em A high fibered power of a
family of 
varieties of general type 
dominates a variety of general type.} Invent. Math. 128 (1997), no. 3,
481--494. 

\bibitem[$\aleph$-K]{A-K} D. Abramovich and K. Karu, {\em  Weak semistable
reduction in characteristic 0}, Invent. math. 139 (2000) 2,
p. 241-273.  



\bibitem[$\aleph$-V1]{fibsurf} D. Abramovich and A. Vistoli, {\em
Complete moduli for fibered surfaces}, 
Recent progress in intersection theory (Bologna, 1997), 1--31, Trends Math.,
Birkh\"auser Boston, 
Boston, MA, 2000.

\bibitem[$\aleph$-V2]{modfam} \ -- \ , {\em
Complete moduli for families over semistable curves}, preprint,
{\tt math.AG/9811059}.

\bibitem[$\aleph$-V3]{stablemaps} \ -- \ , {\em
Compactifying the space of stable maps}, J. Amer. Math. Soc. 15 (2002), no. 1,
27--75. 





\bibitem[Al]{Alexeev:stable-maps}
V. Alexeev, {\em Moduli spaces $M\sb {g,n}(W)$ for surfaces,} in {\it
Higher-dimensional complex varieties (Trento, 
1994)}, 1--22, de Gruyter, Berlin.










\bibitem[D-M]{Deligne-Mumford} 
P. Deligne and  D. Mumford, {\em  The irreducibility of the space of 
curves of given genus,} Inst. Hautes \'Etudes
Sci. Publ. Math. No. {\bf 36} (1969), 75--109. 












\bibitem[J]{DeJong} 
A. J. de Jong, {\em Families of curves and
alterations},  Ann. Inst. Fourier (Grenoble) {\bf 47} (1997), no.~2, 599--621.



\bibitem[Ka]{Karu} K. Karu, {\em Minimal models and boundedness of
stable varieties.} J. Algebraic Geom. 9 (2000), no. 1, 93--109. 






\bibitem[Ko]{Kollar} J.\, Koll\'ar, {\em
Projectivity of complete moduli.} 
J. Differential Geom. 32 (1990), no. 1, 235--268. 


\bibitem[K-SB]{Kollar-Shepherd-Barron} J.\, Koll\'ar and N. Shepherd-Barron, {\em
Threefolds and deformations of surface singularities. } 
Invent. Math. 91 (1988), no. 2, 299--338. 

\bibitem[Ko{\&}al]{flab} J. Koll\'ar et al., { Flips and abundance for
algebraic 
threefolds.   }
Papers from the Second Summer Seminar on Algebraic Geometry held at the
University of Utah, 
Salt Lake City, Utah, August 1991. Ast\'erisque No. 211 (1992). 
Soci\'et\'e Math\'ematique de France, Paris, 1992. pp. 1--258. 




\bibitem[La]{Lanave} G.\ LaNave, {\em Explicit stable models of elliptic
surfaces with sections},  preprint, {\tt math.AG/0205035}








\bibitem[Mo]{Mochizuki} S. Mochizuki, {\it Extending families of
curves over 
log regular schemes.} J. Reine Angew. Math. 511 
(1999), 43--71.







\bibitem[Vi]{Viehweg} E. Viehweg, {\em Quasi-projective moduli for
polarized manifolds.} Ergebnisse der Mathematik und ihrer Grenzgebiete (3),
30. Springer-Verlag, Berlin, 1995. 






\end{thebibliography}
\end{document}